\documentclass[12pt]{article}
\usepackage[sumlimits]{amsmath}
\usepackage{amssymb,amscd}

\usepackage{graphicx}  

\def\nal{\| \hspace{-.10em} |}

\def\1{{\mathbf 1 }}

\catcode`\é=\active\def é{\'e} \catcode`\è=\active\def è{\`e}
\catcode`\ç=\active\def ç{\c{c}} \catcode`\à=\active\def à{\`a}
\catcode`\ê=\active\def ê{\^e} \catcode`\ù=\active\def ù{\`u}
\catcode`\ô=\active\def ô{\^o} \catcode`\û=\active\def û{\^u}
\catcode`\î=\active\def î{\^{\i}} \catcode`\â=\active\def â{\^a}
\catcode`\ö=\active\def ö{\"o}
\catcode`\ë=\active\def ë{\"e}

\def\nnn{{\cal N}}

\def\ph{\varphi}
\def\eps{{\varepsilon}}

\def\Om{{\Omega}}

\def\VaO{V_{\alpha}(\Omega)}

\def \Nal{\mathcal N}

\def \Ual{\mathcal U}
\def \Val{\mathcal V}

\newtheorem{theo}{Theorem}

\newtheorem{prop}[theo]{Proposition}

\newtheorem{exercice}{\bf Exercice}

\def\N{{\mathbb N}}
\def\R{{\mathbb R}}
\def\Z{{\mathbb Z}}

\def\Om{{\Omega}}

\def\ph{\varphi}
\def\eps{{\varepsilon}}

\def\cov{{\rm Cov}}

\def\1{{\mathbf 1}}

\title{ Exponential decay of correlations for a real
 valued dynamical system embedded in $\R^2$}

\author{JAGER Lisette, MAES Jules, NINET Alain}

\begin{document}

\maketitle

%
%

\section*{Abstract }

\noindent
 We study the real valued process  $\{X_t, t\in \N \}$ defined by
$X_{t+2} = \varphi(X_t,X_{t+1})$, where the $X_t$ are bounded.
We aim at proving the decay of correlations for this model, under
regularity assumptions on the transformation $\varphi$.

%
%

\section{Introduction}

\noindent
Since the eighties, the study by statisticians of nonlinear time series has
 allowed to model a great number ot phenomena in Physics, Economics and Finance.
But in the nineties the theory of Chaos became an essential axis of research
for the study of these processes. For an exhaustive review on this subject,
one can consult  Collet-Eckmann \cite{CE} about chaos theory and  Chan-Tong
 \cite{TON1,TON2} about nonlinear time series.
Within this framework, a general model could be written as
$$  X_{t+1}= \ph(X_{t},\dots, X_{t-d+1}) + \varepsilon_t,$$
where $\varphi $ is nonlinear and  $ \varepsilon_t$ is a noise.
We propose a first study of the ``skeleton'' of this model, as Tong calls it, 
beginning with $d=2$ and, more precisely, of the dynamical system induced by
 this model.
Indeed, we consider the model with bounded variables, $X_{t+2}= \ph(X_t,X_{t+1})$,
with $\varphi :\Ual^2 \rightarrow \Ual$ for $\Ual=[-L,L]$ and
 $L \in \R_+^*$, $\varphi$ being defined piecewise on  $\Ual^2$.
This model gives rise to a dynamical system  $ (\Om,\tau ,\mu,T)$ where
 $\mu$ is a measure, invariant under the transformation 
$T :\Om \rightarrow \Om$ and  $\Om$  is a compact subset of $\R^2$.
 Under hypotheses on $\varphi$,
 which imply that $T$ satisfies the hypotheses of Saussol \cite{SAU}, and
if we suppose that $T$ is mixing, we obtain the exponential decay of
 correlations. More precisely, for well-chosen applications $f$ and $h$,
there exist constants  $C=C(f,h)>0$, $ 0<\rho<1$  such that:
 $$ \left|  \int_{\Om} f \circ T^k \,  h \ d\mu - \int_{\Om} f d\mu ~ \int_{\Om}h
 d\mu  \right| \leqslant C \, \rho^k.  $$
This result yields a covariance inequality of the following kind:
$$ \left| \, \cov ( \, f(X_k),h(X_0) \,) \,  \right| \leqslant C \, \rho^k .$$ 
Other ways could certainly be used to get the same result, under different 
hypotheses on the induced system, for
example the method of Young towers  \cite{YOU}. To have a general view on these
different technics, one can read the article of
Alves-Freitas-Luzzato-Vaienti \cite{AFLV}.\\
We finish by giving two examples illustrating our results, a piecewise linear
 one and a nonlinear one.

%
%

\section{Hypotheses and results}\label{hyps-results}

\noindent Let $L \in \R_+^*$. Let   $\varphi : [-L,L]^2 \rightarrow [-L,L]$
\footnote{ To get similar results on $[a,b]$ instead of $[-L,L]$, it suffices
 to conjugate by an affine application} be piecewise defined on  $[-L,L]^2$.
To study the process $\{X_t, t\in \N \}$ defined by  $X_{t+2}= \ph(X_t,X_{t+1})$,
there exist different ways of choosing the induced dynamical system
  $Z_{t+1} = T(Z_t)$ with $Z_t \in \R^2$.
We tried two different approaches, on the one hand the canonical method, setting
 $T(x,y) = (y, \varphi(x,y))$ and on the other hand a double iteration, which
comes down to setting $T(x,y) = (\varphi(x,y),\varphi(y,\varphi(x,y)))$.
The first approach, up to a conjugation, is the most fruitful, the second one
requiring stronger hypotheses and yielding weaker results. We therefore set
$T(x,y) = (\frac{y}{\gamma},\gamma \varphi(x,\frac{y}{\gamma}))$ with
 $Z_t = (X_t,\gamma X_{t+1})$, for a suitable positive $\gamma$.
It then became possible to work in  spaces similar to 
Saussol's $V_{\alpha}$ and to use his results.\\

\noindent 
More precisely, we suppose that the following hypotheses are fulfilled:
\begin{description}
\item[(H1)]
there exists $d \in \N^*$ such that 
$$
[-L,L]^2= \bigcup_{k=1}^d O_k \ \cup \Nal,
$$
where the  $O_k$ are nonempty open sets, $\Nal$ is negligible for the
 Lebesgue measure and the union is disjoint. The edges of the $O_k$ 
can be split into a finite number of smooth components, each one included in a  $C^1$, compact and  one dimensional submanifold of  $\R^2$. 
\item[(H2)]
There exists $\eps_1>0$  such that, for all $k\in \{1,\dots d\}$, there exists
 an application $\ph_k$ defined on $B_{\eps_1}(\overline{O_k}) = \{(x,y) \in \R^2, ~ d((x,y),\overline{O_k}) \leq \eps_1\}$, with values in $\R$, such that  
$\ph_k|_{O_k} = \ph|_{O_k}$.
\item[(H3)]
The application $\ph_k$ is bounded,  $C^{1,\alpha}$ on 
 $B_{\eps_1}(\overline{O_k})$ for a real $\alpha\in ]0,1]$
 \footnote{If $\ph_k$ is $C^2$ on $B_{\eps_1}(\overline{O_k})$, it is
 $C^{1,\alpha}$ on  $B_{\eps_1}(\overline{O_k})$ with $\alpha = 1$}, which means
 that $\ph_k$ is $C^1$ and that there exists  $C_k>0$  such that, for all
 $ (u,v), (u',v')$ in $B_{\eps_1}(\overline{O_k})$, 
$$
\begin{array}{lll}
\displaystyle 
\left|\frac{\partial \ph_k}{\partial u}(u,v) -
\frac{\partial \ph_k}{\partial u}(u',v')\right|
\leq C_k ||(u,v)-(u',v')||^{\alpha}\\
\displaystyle 
\left|\frac{\partial \ph_k}{\partial v}(u,v) -
\frac{\partial \ph_k}{\partial v}(u',v')\right|
\leq C_k ||(u,v)-(u',v')||^{\alpha}.\\
\end{array}
$$
We moreover suppose that there exist  $A>1$ and $M\in ]0,A-1[$  such that   :
$$ 
 \forall (u,v) \in B_{\eps_1}(\overline{O_k}),\qquad
 \left|\frac{\partial \ph_k}{\partial u}(u,v)
\right|\geq A,\quad
\left|\frac{\partial \ph_k}{\partial v}(u,v)
\right| \leq M,
$$
to ensure the expansion.
\item[(H4)]
The open sets $O_k$  satisfy the following geometrical condition:
 \footnote{In suitable cases, this hypothese can be replaced by a weaker but 
simpler one : for all  points $(u,v)$ and $(u',v)$ in
 $B_{\eps_1}(\overline{O_k})$, the segment $[(u,v),(u',v)]$ is included in 
$B_{\eps_1}(\overline{O_k})$}
 For all  $(u,v)$ and $(u',v)$ in $ B_{\eps_1}(\overline{O_k})$,  there exists 
 a $C^1$  path  $\Gamma=(\Gamma_1,\Gamma_2) : [0,1]\rightarrow  B_{\eps_1}(\overline{O_k})$ $C^1$ joining $(u,v)$ and $(u',v)$, whose gradient does not
 vanish and which satisfies
\begin{equation*}
\forall t\in ]0,1[, \left|\Gamma_1'(t)\right|
 >\frac{M}{A}  \left|\Gamma_2'(t)\right| .
\end{equation*}
\item[(H5)]
Let  $Y \in \N^*$ be  the maximal number of $C^1$ components of 
 $\Nal$ meeting at one point and  set
$$
s = \left( \frac{  2A+M^2 -M\sqrt{M^2+4A} }{2} \right)^{-1/2} < 1.
$$
One supposes that 
$$
\eta:=  s^{\alpha} + \frac{8 s}{\pi(1-s)}Y <1.
$$
\end{description}

\noindent We set $\displaystyle \gamma = \frac{1}{\sqrt{A}}<1$ and,  for all
 $k \in \{1, ..., d\}$, we denote by  $U_k$ (resp. $W_k$, $\Nal'$) 
 the image of $O_k$ 
(resp. $B_{\eps_1}(\overline{O_k})$, $\Nal$) under the compression which
associates  $(u,\gamma v)$ with each $(u,v) \in \R^2$.

 The set  $\Om = [-L,L] \times [-\gamma L, \gamma L]$, on which we shall be working, is the image of  $[-L,L]^2$ under the same compression.
\\
For every non negligible Borel set $S$ of  $\R^2$, for every 
 $f \in L^1_m(\R^2,\R)$,  set
$$ Osc(f,S) = \underset{S}{Esup} f - \underset{S}{Einf} f, $$
where $\underset{S}{Esup}$ and $\underset{S}{Einf}$ are the essential
supremum and infimum with respect to the Lebesgue measure $m$.
One then defines: 
$$
|f|_{\alpha}= \sup_{0<\eps<\eps_1}\eps^{-\alpha} \int_{\R^2}  {\rm Osc}(f,B_{\eps}(x,y))\ dxdy \qquad,  \qquad \| f \|_{\alpha} = \| f \|_{L^1_m} + |f|_{\alpha}
$$
and the set
 $V_{\alpha} = \{ f \in L^1_m(\R^2,\R), ~ \| f \|_{\alpha} < +\infty \}$.\\

\noindent 
Let us introduce similar notions on
 $\Om$ :  for every $0 < \eps_0 < \gamma \eps_1$,  for every
 $g \in L^{\infty}_m(\Om,\R)$, one defines
$$
N(g,\alpha,L) = \sup_{0<\eps<\eps_0}\eps^{-\alpha} \int_{\Om} {\rm Osc}(g,B_{\eps}(x,y) \cap \Om)\ dxdy.
$$
One then sets:
$$
||g||_{\alpha,L}=  N(g,\alpha,L) + 16 (1+\gamma) \eps_0^{1-\alpha} L ||g||_{\infty} + ||g||_{L^1_m}.
$$
The function $g$ is said to belong to $\VaO$ if the above expression is
 finite. The set  $\VaO$ does not depend on the choice of
 $\eps_0$, whereas  $N$ and $\| . \|_{\alpha,L}$ do.\\
 There exist relationships between these two sets. Indeed, thanks to
 Proposition 3.4 of  \cite{SAU}, one can prove the
 following result:

\begin{prop}
\begin{enumerate}
\item
If $g \in \VaO$ and if one extends $g$ as a function denoted by $f$, setting 
$f(x,y)=0$ if $(x,y)\notin \Om$, then $f \in V_{\alpha}$ and
$$
\| f \|_{\alpha} \leq \| g \|_{\alpha,L}.
$$
\item
Let $f$ be in $ V_{\alpha}$. Set $g= f \1_{\Om}$. Then $g\in \VaO$ and one has
$$
 \| g \|_{\alpha,L}\leq \left( 1+16 (1+\gamma) L
 \frac{\max(1,\eps_0^{\alpha})}{\pi \eps_0^{1+\alpha}}\right)
\| f \|_{\alpha}.
$$
\end{enumerate}
\end{prop}

\noindent 
Under the above hypotheses  (H1) to (H5), one obtains a first result:
\begin{theo}\label{resconjug}
Let $T$ be the transformation defined on $\Om$ by : $\forall  (x,y) \in U_k$ :
$$
T(x,y) = T_k(x,y)= \left( \frac{y}{\gamma},\
 \gamma \ph_k(x, \frac{y}{\gamma})\right).
$$
Keeping the same formula, one extends the definition of  $T_k$  to $W_k$. Then
\begin{enumerate}
\item 
The  Frobenius-Perron operator $P : L^1_m(\Om) \rightarrow L^1_m(\Om)$ 
associated with $T$ has a finite number of eigenvalues 
 $\lambda_1,\dots,\lambda_r$ of modulus one.
\item 
For each $i\in\{1,\dots,r\}$, the eigenspace
 $E_i=\{ f\in L^1_m(\Om) \ : \ Pf=\lambda_i f\}$ associated with the eigenvalue
  $\lambda_i$ is finite dimensional and included in  $V_{\alpha}(\Om)$.
\item 
The operator $P$ decomposes as
$$
P= \sum_{i=1}^r \lambda_i P_i + Q,
$$
where the $P_i$ are projections on the spaces $E_i$, $\nal P_i \nal_1\leq 1$
and $Q$ is a linear operator defined on $ L^1_m(\Om)$, satisfying 
 $Q(\VaO) \subset \VaO$, $\sup\limits_{n \in \N^*}\nal Q^n \nal_1<\infty$ and
  $\nal Q^n \nal_{\alpha,L} = O(q^n)$ when $n \rightarrow +\infty$
for an exponent $q\in ]0,1[$. Moreover, $P_iP_j=0$ if $i\neq j$, 
$P_iQ=QP_i=0$  for all $i$.
\item 
The number $1$ is an eigenvalue of $P$. Set  $\lambda_1=1$, let 
 $h_*=P_1 \1_{\Om}$ and let $d\mu=h_* ~ dm$.
Then $\mu$ is the greatest absolutely continuous invariant measure (ACIM) of
 $ T$, that is to say:
if $\nu <<m$ and if $\nu$ is $T$-invariant, then $\nu<<\mu$.
\item 
The support of $\mu$ can be decomposed into a finite number of disjoint
measurable sets, on which a power of $T$ is mixing. More precisely  for all
 $j \in \{1,2,\dots, \dim(E_1)\}$, there exist an integer $L_j \in \N^*$ and 
 $L_j$  disjoint sets $W_{j,l}$ $(0 \leq l \leq L_{j}-1)$ satisfying 
 $T(W_{j,l})=W_{j,l+1 \mod(L_j)}$ and $T^{L_j}$ is mixing on every  $W_{j,l}$.
We denote by $\mu_{j,l}$ the normalized restriction of  $\mu$ to $W_{j,l}$,
defined by
$$
\mu_{j,l}(B)= \frac{\mu(B\cap W_{j,l})}{\mu(W_{j,l})}, \ d\mu_{j,l} =
 \frac{h^* \1_{W_{j,l}} }{\mu(W_{j,l})} dm.
$$
The fact that  $T^{L_j}$ is mixing on every $W_{j,l}$ means that, for all
 $f \in L^1_{\mu_{j,l}}(W_{j,l})$ and all $h \in L^{\infty}_{\mu_{j,l}}(W_{j,l})$,
$$ \lim\limits_{t \rightarrow + \infty} <T^{tL_j} f,h>_{\mu_{j,l}} = 
<f,1>_{\mu_{j,l}} <1,h>_{\mu_{j,l}} $$  with the notations (indifferently employed)
 $<f,g>_{\mu'} = \mu'(fg) = \int f g ~ d\mu'$. 
\item Moreover, there exist  $C>0$ and $0<\rho < 1$  such that, for all
 $h$ in $\VaO$ and $f\in L^1_{\mu}(\Om)$, one has
$$
\left| \int_{\Om} f \circ T^{ k \times ppcm(L_i)}  h \ d\mu -\sum_{j=1}^{\dim(E_1)} \sum_{l=0}^{L_j-1} \mu(W_{j,l}) <f,1>_{\mu_{j,l}} <1,h>_{\mu_{j,l}} \right| \leq C  ||h||_{\alpha,\Om} ||f||_{L^1_{\mu}(\Om)}~\rho^{k}.
$$
\item If, moreover, $T$ is mixing\footnote{ which is equivalent to: if $ 1$ is
 the only eigenvalue of $P$ with modulus one and if it is simple}, then the preceding result can be written as follows:
  there exist  $C>0$ and $0<\rho < 1$  such that, for all $h$ in $\VaO$
 and $f\in L^1_{\mu}(\Om)$, one has:
$$ \left|  \int_{\Om} f \circ T^k \, h \ d\mu - \int_{\Om} f d\mu ~ \int_{\Om} h
 d\mu  \right|  \leq C  ||h||_{\alpha,\Om}~||f||_{L^1_{\mu}(\Om)}~\rho^{k}. $$
\end{enumerate}
\end{theo}

\noindent Now let us come back to the initial problem and try to deduce
from this result the invariant law associated with $X_t$. If 
 $(X_t)_t$ is defined by  $X_0, X_1$ (valued in $[-L,L]$)  and the recurrence
relation $X_{t+2}= \ph(X_t,X_{t+1})$, one sets $Z_t=(X_t,\gamma X_{t+1})$. 
Then $(Z_t)_t$ satisfies the recurrence relation  $Z_{t+1}=T(Z_t)$,
which implies the following result:

\begin{theo}
Suppose that the random variable  $Z_0=(X_0,\gamma X_{1})$ has the density $h_*$.
Then  $Z_t$ has the density $h_*$ and for all $t \in \N$, $X_t$ has the density
\begin{equation}
f : x\mapsto \int_{[-\gamma L,\gamma L]} h_*(x,v)\ dv
= \gamma \int_{[- L, L]} h_*(u,\gamma x)\ du. \label{1}
\end{equation}
\end{theo}

\noindent 
Indeed, since $Z_t= (X_t,\gamma X_{t+1})$ has the density $h_*$, one proves that
 $X_t$ has the density $f$ by computing the first marginal distribution.
Computing the second one yields that $\gamma X_{t+1}$ has the density $g=g(y)$
 defined by
$$
g(y)= \int_{[-L,L]} h_*(u,y)\ du.
$$
This implies that $ X_{t+1}$ has the density $y\mapsto \gamma g(\gamma y)$.
 But $Z_{t+1}$ has the density $h_*$ as well. Therefore $X_{t+1}$
has the density given by the first marginal distribution, which proves 
 the equality (\ref{1}).
\\

\noindent If  $F$ is defined on $[-L,L]$, we denote by $Tr ~ F$ 
the function defined, on $\Omega$, by $Tr ~ F(x,y)= F(x)$.\\
\noindent 
One then obtains the following result, which is a direct consequence of
the sixth point of
 Theorem \ref{resconjug}, applied to $Tr ~ F$ and $Tr ~ H$ :

\begin{theo}\label{th4}
For every Borel set $B$ and every interval  $I$, if $(X_0,X_1)$ has the 
invariant distribution, then
$$
\begin{array}{c}
\displaystyle \left\vert P\left( X_{ k \times ppcm(L_i)}\in B, X_0 \in I \right)  -\sum_{j=1}^{\dim(E_1)} \sum_{l=0}^{L_j-1} \mu(W_{j,l}) <Tr ~ \1_B,1>_{\mu_{j,l}} <1,Tr ~ \1_I>_{\mu_{j,l}} \right| \\
\\
\leq 16 (1+\gamma) ~ C ~ L^3 ~ (10 \eps_0^{1-\alpha} + L) ~ \rho^{k }.
\end{array}
$$
More generally, let $F$, defined and measurable on $[-L,L]$, be  such that 
  $Tr ~ F$ belongs to  $L^1_{\mu}(\Om)$. Let $H \in L^{\infty}_m([-L,L])$ be
  such that   $\displaystyle \sup_{0<\eps<\eps_0}\eps^{-\alpha} \int_{[-L,L]}
 {\rm Osc}(H,]x-\eps,x+\eps[ \cap [-L,L])\ dx < +\infty$.\\
 Then $Tr ~ H \in V_{\alpha}(\Om)$ and 
$$
\left\vert E( F( X_{k \times ppcm(L_i)}) H(X_{0}))
 -\sum_{j=1}^{\dim(E_1)} \sum_{l=0}^{L_j-1} \mu(W_{j,l})
 \mu_{j,l}(Tr ~ F)\mu_{j,l}(Tr ~ H) \right| 
\leq C(F,H) ~ \rho^k $$
with
$$ \begin{array}{l}
\displaystyle
C(F,H) = C ~ L ~ ||Tr ~ F||_{L^1_{\mu}}\left(2\gamma \sup_{0<\eps<\eps_0}\eps^{-\alpha} \int_{[-L,L]} {\rm Osc}(H,]x-\eps,x+\eps[ \cap [-L,L])\ dx \right. 
\\
\displaystyle + 16 (1+\gamma) ~ \eps_0^{1-\alpha} ||H||_{L^{\infty}_m([-L,L])} + 2\gamma ~ ||H||_{L^1_m([-L,L])}\bigg).
\end{array}
$$
If, moreover, $T$ is mixing, then:
$$ |Cov(F(X_k),H(X_0))| \leq C(F,H) ~ \rho^k. $$
\end{theo}

%
%

\section{Proofs}\label{demos}

\noindent Theorem \ref{resconjug}  is a consequence of Theorems 5.1 and 6.1 of
 \cite{SAU}.  The difficulty is proving that $T$ satisfies Hypotheses
  (PE1) to (PE5).\\

\noindent To check that (PE2) is satisfied, we first prove that $T_k$ is a
$C^1$ diffeomorphism from  $W_k$ on $T_k(W_k)$. Hypothesis (H3) about
 $\displaystyle \frac{\partial \ph_k}{\partial u}$ ensures that $T_k$ is a
 local diffeomorphism. To establish the injectivity, let us consider two
different points $(x,y)$ and $(x',y')$ of $W_k$, whose image under $T$
is the same. One then has $y=y'$ and $\ph_k(x',y/\gamma)= \ph_k(x,y/\gamma)$.
Using the geometrical hypothesis (H4) and applying the Mean Value Theorem to 
the application
 $t \mapsto \ph_k(\Gamma_1(t),\Gamma_2(t) )$, one obtains a contradiction.

\noindent The regularity hypotheses on $\ph_k$ (and consequently on $T_k$)
imply that $\det(DT_k^{-1})$ is Hölder contiuous for the exponent $\alpha$, 
on a  suitably restricted domain. One can prove that there exist, for every
$k$,  a real number $\beta_{k}>0$, an open set $\Val_k$ with compact closure
and a constant  $c_{k}$  such that  
\begin{itemize}
\item $\overline{U_k}\subset  \Val_{k} \subset \overline{\Val_{k}} \subset W_k$ ;
\item $B_{\beta_{k}}(T_k(U_k)) \subset T_k(\Val_{k}) ;$
\item
for every $\eps<\beta_{k}$, every $z\in T_k(\Val_{k})$ and all
 $x,y\in B_{\eps}(z)\cap T_k(\Val_{k})$, one has
$$
\Big| \det(DT_k^{-1}(x))- \det(DT_k^{-1}(y)) \Big| \leq c_{k} \Big| \det(DT_k^{-1}(z)) \Big|\eps^{\alpha}.
$$
\end{itemize}

\noindent Setting $\beta = \min\limits_k \beta_{k} >0$ and
 $c = \max\limits_k c_{k} >0$, one gets constants which are valid for all
 $k \in \{1, \hdots , d\}$. Hence  (PE2) is satisfied.
\\

\noindent This allows to fix the open set with which we are going to work:
  there exists  $\eps_2>0$  such that, for all $k\in\{1,\dots,d\}$,
$B_{2\eps_2}(\overline{U_k}) \subset \Val_{k} \subset W_k$.
From now on,  $V_k=B_{\eps_2}(\overline{U_k})$. The set $T_k(V_k)$ is open and 
$T_k(\overline{U_k})$ is a compact set included in $T_k(V_k)$. One can find
 $\eps_0^1 >0$  such that $B_{\eps_0^1}( T_k(\overline{U_k})) \subset T_k(V_k)$
  for all $k$. Hypothesis (PE1) is thus verified.\\
\\
Hypothesis (PE3) is clearly satisfied because 
 $\displaystyle \Om= \bigcup_{k=1}^d  U_k\ \cup \Nal'$ is the disjoint union of
open sets and of a negligible set.
\\

One treats (PE4) in two steps  :
first one proves an expansion result, in the case when the arguments 
in   $\Val_k$ are near (Proposition \ref{dilatance1}), then one proves 
(PE4) itself, which is an expansion result in the case when the images (in 
 $T_k(V_k)$ ) are near.

\begin{prop}\label{dilatance1}
Let $(x,y)$ and $(x',y')\in \Val_{k}$  be such that the
 segment $[(x,y),(x',y')]$ is included in $\Val_{k}$. Then
$$
|| T_k(x,y)-T_k(x',y')||^2 \geq \frac{1}{s^2} || (x,y)-(x',y')||^2.
$$
\end{prop}

\noindent{\it Proof:}
Applying the Mean Value Theorem to the application defined on  $[0,1]$ by
 $t \mapsto \ph_k(x+t(x'-x), \frac{1}{\gamma} (y+t(y'-y))$ gives a number
  $c\in ]0,1[$  such that 
$$
 || T_k(x,y)-T_k(x',y')||^2 = (x'-x,y'-y) B \left( \begin{array}{lll}
x'-x\\
y'-y
\end{array}\right)
$$
where
$$
B = \left( \begin{array}{lll}
\displaystyle
\gamma^2 
\left(\frac{\partial \ph_k}{\partial u}(x_c,\frac{1}{\gamma}y_c)\right)^2
&
\displaystyle
\gamma \frac{\partial \ph_k}{\partial u}(x_c,\frac{1}{\gamma}y_c)\frac{\partial \ph_k}{\partial v}(x_c,\frac{1}{\gamma}y_c)\\
\displaystyle
\gamma \frac{\partial \ph_k}{\partial u}(x_c,\frac{1}{\gamma}y_c)\frac{\partial \ph_k}{\partial v}(x_c,\frac{1}{\gamma}y_c)  &
\displaystyle \frac{1}{\gamma^2} +
\left(\frac{\partial \ph_k}{\partial v}(x_c,\frac{1}{\gamma}y_c)\right)^2
\end{array}\right)
$$
with $(x_c,y_c) = (x+c(x'-x),y+c(y'-y))$.\\
The matrix $B$ is real and symmetrical. Set 
$$
\begin{array}{lll}
\displaystyle
\xi_1 & = {\rm Tr}(B) & \displaystyle = 
 \frac{1}{\gamma^2} +
\left(\frac{\partial \ph_k}{\partial v}(x_c,\frac{1}{\gamma}y_c)\right)^2
+ \gamma^2 
\left(\frac{\partial \ph_k}{\partial u}(x_c,\frac{1}{\gamma}y_c)\right)^2 \\
\xi_2 & =\det(B) &\displaystyle = \left(\frac{\partial \ph_k}{\partial u}
(x_c,\frac{1}{\gamma}y_c)\right)^2.
\end{array}
$$
 We now prove that the eigenvalues of $B$ are greater than $\frac{1}{s^2}$.
Indeed, the map $\zeta :\R^2 \rightarrow \R^2$ defined by
 $\zeta(x,y)= (x+y,xy)$ is a bijection from
$$
V''_s=\{ (x,y) \in \R^2 : s^{-2} \leq x\leq y\}
$$
to
$$
 \zeta(V''_s)= \{ (\xi_1,\xi_2) \in (\R^*_+)^2 : \xi_1 \geq 2s^{-2},\ 
\xi_2\geq s^{-2}(\xi_1-s^{-2}),\ \xi_2 \leq \frac{\xi_1^2}{4} \}.
$$
One just has to check that  $(\xi_1,\xi_2)$ is in $\zeta(V''_s)$ to obtain
 the result.\\
Now since $B$ has real eigenvalues, the discriminant of its characteristic
 polynomial is nonnegative. Consequently  $4\xi_2\leq \xi_1^2$. The conditions
 on  $A$ and $M$ and the choice of $s$ and $\gamma$ imply that the other
 inequalities are satisfied. 
\\ It follows that eigenvalues of the matrix $B$ are greater than or equal to
 $s^{-2}$. Hence $ || T_k(x,y)-T_k(x',y')||^2 \geq \frac{1}{s^2} || (x,y)-x',y')||^2$, which completes the first step.  \hfill $\square$\\\\

\noindent
Compacity arguments prove that there exists  $\eps_0^2 > 0$  such that, 
 for all $(x,y)\in \overline{V_k}$,
$$
B_{\eps_0^2}(T_k(x,y)) \subset T_k(B_{\eps_2}(x,y)).
$$

\begin{prop}\label{dilatance2}
Set $\eps_0= \min(\eps_0^1,\eps_0^2) > 0$. Recall that $\overline{U_k}
 \subset V_k \subset \overline{V_k} \subset \Val_k \subset W_k$. Then:
\begin{itemize}
\item
For all $(u_1,v_1), (u_2,v_2) \in T_k(V_k)$ satisfying 
 $d((u_1,v_1), (u_2,v_2))<\eps_0$, the following inequality holds:
$$
s^2 ~ d((u_1,v_1), (u_2,v_2)) >
 d(T_k^{-1}(u_1,v_1),(T_k^{-1} (u_2,v_2)),
$$
with $\displaystyle s^2 < 1$.
\item
$B_{\eps_0}( T_k(\overline{U_k})) \subset T_k(V_k)$.
\end{itemize}
\end{prop}

\noindent{\it Proof :}
The second assertion comes from the fact that $\eps_0\leq \eps_0^1$  and from
what we have obtained in  (PE1). \\
The first assertion implies Condition (PE4) of Saussol. To prove it, let
 $(u_1,v_1), (u_2,v_2) \in T_k(V_k)$ satisfy $d((u_1,v_1), (u_2,v_2))<\eps_0$.
 Let $(x,y)=T_k^{-1}(u_1,v_1)$ be in $ V_k$. According to the preceding remark,
as  $\eps_0$ is smaller than $\eps_0^2$,
$$
(u_2,v_2) \in B_{\eps_0}(T_k(x,y)) \subset T_k(B_{\eps_2}(x,y)).
$$
Hence $(x',y')=T_k^{-1}(u_2,v_2)\in B_{\eps_2}(x,y)\subset \Val_k$.
According to the Proposition  \ref{dilatance1}, 
$$
d((u_1,v_1),(u_2,v_2))^2 = 
|| T_k(x,y)-T_k(x',y')||^2 > \sigma || (x,y)-(x',y')||^2,
$$
which proves the result. \hfill $\square$\\

\noindent To conclude, Hypothesis  (PE5) is a consequence of Lemma  2.1 of
 Saussol and of Hypothesis  (H5).\\

\noindent Since the hypotheses (PE1) to (PE5) are verified, Theorem 5.1 of
 \cite{SAU} implies the properties 1 to 5 of Theorem \ref{resconjug} about 
 $V_{\alpha}$ and $L^1_m$. But, if $f \in E_i$, $f$ is equal to $0$ on
 $\Omega^c$ and then $f$ belongs to $L^1_m(\Omega)$ and to $V_{\alpha}(\Omega)$.\\

\noindent To prove the point 6, we apply Theorem 6.1 of \cite{SAU} on every
$W_{j,l}$, on which a power of $T$ is mixing. 
Using the notations of Point 5 of Theorem 5.2 of \cite{SAU}, one obtains the
 existence of real constants $C>0$ and $\rho\in ]0,1[$ such that,  for all
 $(j,l)$ satisfying $ 1\leq j\leq \dim(E_1)$, $0\leq l\leq L_j -1$, for every 
function $f\in L^1_{\mu_{j,l}}(\Om)$ and for every function  $h\in V_{\alpha}(\Om)$,
$$
\left\vert \int_{\Om}(f- \mu_{j,l}(f))\circ T^{k L_j} h ~ d\mu_{j,l}
\right\vert \leq C ||f-\mu_{j,l}(f))||_{L^1_{\mu_{j,l}}} ||h||_{\alpha,L}\rho^{k } .
$$
Let then  $h$ be in $  V_{\alpha}(\Om)$ and  $f$ be in $ L^1_{\mu}(\Om)$ 
(so that $f\in L^1_{\mu_{j,l}}(\Om)$ for every  $j,l$). Taking the smallest
 common multiple $L'$ of the  $L_j$ and summing the above inequalities, with
$k$ replaced with  $k \frac{L'}{L_j}$, one gets
$$
\left| \int_{\Om} f \circ T^{ k L'}  h \ d\mu -\sum_{j=1}^{\dim(E_1)}
 \sum_{l=0}^{L_j-1} \mu(W_{j,l}) \mu_{j,l}(f)\mu_{j,l}(h) \right| 
\leq C  ||h||_{\alpha,\Om}||f||_{L^1_{\mu}} \rho^{k } .
$$

\noindent Point 7 is a direct consequence of Point 6, since
 $\dim(E_1) = 1$ and $L_1=1$. This completes the proof of Theorem 
\ref{resconjug}.\hfill $\square$ \\

\noindent Now let us prove Theorem \ref{th4}.
If $\left( \begin{array}{lll} X_{0}\\ \gamma X_{1} \end{array}\right)$ has the
 distribution $\mu$, then so does $\left( \begin{array}{lll} X_{k}\\ \gamma X_{k+1} \end{array}\right)$. If $f\in L^1_{\mu}(\Om)$ and if $h\in V_{\alpha}(\Omega)$, we then have:
$$
\left\vert E\left( f \left( \begin{array}{lll} X_{kL'}\\ \gamma X_{kL'+1} \end{array}\right) h\left( \begin{array}{lll} X_{0}\\ \gamma X_{1} \end{array}\right)\right)
 -\sum_{j=1}^{\dim(E_1)} \sum_{l=0}^{L_j-1} \mu(W_{j,l}) \mu_{j,l}(f)\mu_{j,l}(h) \right| \leq C ||f||_{L^1_{\mu}} ||h||_{\alpha,\Om}\rho^{k } .
$$
\bigskip

\noindent  In order for  $Tr~ H$ to belong to $V_{\alpha}(\Omega)$, it is 
sufficient and necessary that $H$ belongs to $L^{\infty}([-L,L],m)$ and satisfies
$$
 \sup_{0<\eps<\eps_0}\eps^{-\alpha} 
\int_{[-L,L]} {\rm Osc}(H,]x-\eps,x+\eps[\cap [-L,L])\ dx  <\infty.
$$
Moreover,
$$
\begin{array}{lll}
|| Tr ~ H||_{\alpha,\Om} 
 & \displaystyle =
2\gamma L \sup_{0<\eps<\eps_0}\eps^{-\alpha} \int_{[-L,L]}
 {\rm Osc}(H,]x-\eps,x+\eps[ \cap [-L,L])\ dx\\
& \displaystyle
 + 16(1+\gamma) L \eps_0^{1-\alpha} ||H||_{L^{\infty}_m([-L,L])}
 + 2\gamma L ||H||_{L^1_m([-L,L])}.
\end{array}
$$
Thus if $H$ satisfies these conditions and if $F$ is such that $Tr~F$ belongs
 to $L^1_{\mu}(\Om)$, for example if $F$ is measurable and bounded on 
$[-L,L]$, one has

$$
\begin{array}{lllll}
\displaystyle
\left\vert E( F( X_{k \times L'}) H(X_{0}))
 -\sum_{j=1}^{\dim(E_1)} \sum_{l=0}^{L_j-1} \mu(W_{j,l})
 \mu_{j,l}(Tr ~ F)\mu_{j,l}(Tr ~ H) \right| 
\\
\displaystyle\leq C 
 ||Tr ~ F||_{L^1_{\mu}}\left(
2\gamma L \sup_{0<\eps<\eps_0}\eps^{-\alpha} \int_{[-L,L]}
 {\rm Osc}(H,]x-\eps,x+\eps[ \cap [-L,L])\ dx \right. 
\\
\displaystyle + 16 (1+\gamma) L \eps_0^{1-\alpha} ||H||_{L^{\infty}_m([-L,L])} + 2\gamma L ||H||_{L^1_m([-L,L])}\bigg)\rho^{k } .
\end{array}
$$

\noindent In particular, if $H$ is the characteristic function  of an interval
 and if $F$  is the characteristic function  of a Borel set, we obtain the
first assertion of Theorem \ref{th4}.

%
%

\section{Examples}\label{exemples}

\subsection{A nonlinear example}

\noindent For all $k \in \Z$ we denote by $f_k$ the polynomial function
 $\displaystyle f_k(x) = - \frac{71}{2} x^2 - 214 x + k - \frac{1}{2}$.\\
For all $ -179 \leq k \leq 250$, one defines the open set $O_k$ by :
$$ O_k = \{ (u,v) \in ]-1,1[^2 ~ / ~ f_k(u) < v < f_{k+1} (u) \}. $$

\noindent  We consider the applications $\ph_k$  defined on
 $B_1(\overline{O_k})$ ($\eps_1 = 1$)  for all $-179 \leq k \leq 250$ by :
$$ \ph_k(u,v) = 2v - 2 f_k(u) - 1. $$
One defines $\ph : [-1,1]^2 \rightarrow [-1,1]$ almost everywhere by setting
 $\ph_{|_{O_k}} = \ph_{k|_{O_k}}$  for all $-179 \leq k \leq 250$.
We now make sure that these functions and open sets satisfy the conditions
specified in Section \ref{hyps-results}.

\noindent
The condition about the open sets is easily verified, since 
$[-1,1]^2 \backslash \bigcup\limits_{k=-179}^{250} O_k $ is the union of 
segments and parabolic arcs. Moreover, the maximal number of arcs meeting at
 one point is  $Y=3$. \\
The regularity conditions are satisfied, because the  $\ph_k$ are smooth on 
$B_1(\overline{O_k})$. Set  $\alpha = 1$. The partial derivatives satisfy
the following inequalities:  for all $-179 \leq k \leq 250$ and  all
 $(u,v) \in B_1(\overline{O_k})$ one has:
$$ \left| \frac{\partial \ph_k}{\partial v} (u,v) \right| = 2 = M $$
and
$$ \left| \frac{\partial \ph_k}{\partial u} (u,v) \right| = 2 | 71 u + 214| > 2 (214 - 71 (1 + 1)) = 144 = A > M+1. $$

\noindent In this case, $\gamma = \frac{1}{12}$. A computation shows that 
 $s \leq \frac{1}{10}$ and $\eta < 1$.\\
One sets $\Om= [-1,1]\times\left[ - \frac{1}{12},\frac{1}{12}\right]$ and for
 all $-179 \leq k \leq 250$ one defines the open sets
$$
U_k=  \{(x,y)\in \mathring{\Om} : f_k(x) < 12 y < f_{k+1}(x) \}.
$$
We obtain the applications
$$
T_k(x,y)= (12 y, \frac{2}{12}(12 y -f_k(x))- \frac{1}{12}).
$$

\noindent If $-177 \leq k\leq 248$, $T_k(U_k)= \mathring{\Om}$ and $T_k$ 
is a bijection from $U_k$ on $\mathring{\Om}$.\\
Otherwise, one can check that $T_{-178}$ is a bijection from $U_{-178}$ on
 $\displaystyle \Om_1 \cup \Om_2$, where $\Omega_1$ is the open subset of
 $\mathring{\Omega}$ above the  line having the equation $y= \frac{2x+1}{12}$, 
$\Omega_3$, the open subset under the line
having the equation $y=\frac{2x-1}{12}$ and
 $\Omega_2$, the one between both lines. One has similar relations for
 $k=-177$, $249$ and $250$ and other subsets of $\Om$.\\

\noindent Finally, the simple version of the geometrical condition is satisfied
 (the open set contains the horizontal segment).\\
The transformation $T$ therefore admits an invariant density $h_*$.\\

\noindent Let $P$ be the Frobenius-Perrron operator associated with $T$.
One can prove that the constant functions are not invariant by $P$ and
consequently that $h_*$ is not constant. Indeed, set
$$
\psi_k(x,y)= (214)^2-71(2x-12y)+142 k.
$$
Then $Ph$ can be written as $\displaystyle Ph(x,y)= \sum_{k=a}^{b} h(T_k^{-1}(x,y)) \frac{1}{2 \sqrt{\Psi_k(x,y)}}$, with $(a,b)=(-179,248)$ if $(x,y) \in \Omega_1$, $(a,b) = (-178,249)$ si $(x,y) \in \Omega_2$, $(a,b) = (-177,250)$ if $(x,y) \in \Omega_3$.\\

\noindent We now verify that $P1\neq 1$. Suppose that $h=1$ and set $z=x-6y$. \\
If $(x,y)\in \Om_3$, $z \in ]-\frac{3}{2}, -\frac{1}{2}[$. The function
 $\displaystyle z\mapsto \sqrt{ (214)^2-142 z +142 k}$ is strictly decreasing
 on $]-\frac{3}{2}, -\frac{1}{2}[$. Therefore 
$\displaystyle z\mapsto  \frac{1}{2 \sqrt{ (214)^2-71(2x-12y)+142 k}}$
is strictly increasing on $]-\frac{3}{2}, -\frac{1}{2}[$ and $P1$ is not 
constant.

\subsection{A piecewise linear example}

\noindent This example can be useful to create a generator of  pseudo random
 numbers in $[-L,L]$.\\
In this section, $a$ and $b$ are positive or negative integers, $L$ is a positive integer or half integer. \\
One denotes by $\Ual^2$ the square $[-L,L]^2$. For all $n\in \Z$, the open set
 $\Omega_n$ is defined by
$$
\Omega_n=\{ (u,v)\in ]-L,L[^2 \ : \  av+bu \in ] (2n-1)L,(2n+1)L[\}.
$$
One denotes by $\Delta_n$  the line having the equation $ av+bu=(2n-1)L$.
One defines  $\varphi_n$ on $\R^2$ by
$$
\varphi_n(u,v)=  av+bu-2nL.
$$
Then $\left.\varphi_n\right|_{\Omega_n}$ is valued in $]-L,L[$ and we set
$$
\forall (u,v)\in \Omega_n,\ \varphi(u,v) = \varphi_n(u,v).
$$
We impose the following condition, with
 $\displaystyle S= 1+ \frac{48}{\pi} + \frac{288}{\pi^2} + 
\frac{4}{\pi} \left(1+\frac{12}{\pi}  \right)\sqrt{6\pi +36}$,
$$ |a|< \frac{ |b| -S}{\sqrt{S}}. $$
\vskip 0.5 cm

\noindent One verifies that the conditions of
Section  \ref{hyps-results} are fulfilled. \\
The square $\Ual^2$ is the disjoint union of a finite number of open sets
 $\Om_n$ and of a negligible set composed of a finite number of
segments.\\
The maximal number of these segments meeting at one point is  $Y=3$.\\
For every $\eta>0$, the open sets $ B_{\eta}(\Omega_n)$ are convex, hence
 contain the horizontal segment joining two points having the same ordinate
and the geometrical condition is satisfied.\\

The applications $\varphi_n$  are smooth on $ B_{\eta}(\Omega_n)$. We set
 $\alpha=1$. Moreover $\varphi_n(\Omega_n)\subset [-L,L]$.\\
We set $M=|a|$, $A=|b|$, so that the partial derivatives satisfy the required
inequalities. The upper bound of $|a|$ implies that $0<M<A-1$.\\
One sets $\gamma= |b|^{-1/2}$ (it is the compression coefficient) and one 
checks that   $\eta<1$.\\

\noindent 
We determine for which integers $n$ the line  $\Delta_n$ crosses the square.
One can see that  $\Delta_n\cap \Ual^2 \neq \emptyset$ if and only if 
$$
\frac{-|a|- |b|+1}{2}\leq n\leq \frac{|a|+|b|+1}{2}.
$$

\noindent One defines $\nnn(a,b)$ as the set of indices $n$  such that
a nonempty $\Omega_n$  intersects $\Ual^2$.

\noindent Set
$$
\Omega= [-L,L]\times [-\gamma L,\gamma L]
$$
and
$$
\Omega_{n,a}=
\{ (x,y)\in \mathring{\Omega}\ :\ a \sqrt{|b|} y + b x -2nL \in ]-L,L[ \}.
$$
Set
$$
\begin{array}{lllll}
T_n & : & \Omega_{n,a} & \rightarrow & \Omega\\
&& (x,y) & \mapsto &  ( \sqrt{|b|} y, ay+ \frac{b}{ \sqrt{|b|} }x
- \frac{2nL}{ \sqrt{|b|} } ). \\
\end{array}
$$
One defines $T$ almost everywhere from $\Omega$ in $\Omega$ by setting
 $\left. T \right|_{\Omega_{n,a}} = T_n$.\\

\noindent The Frobenius-Perron operator $P$ associated with $T$ has 
the following expression: 
$$
Ph(x,y)=\frac{1}{|b|}
\sum_{n\in \nnn(a,b)} \1_{(x,y)\in T_n(\Omega_{n,a}) } h(T_n^{-1}(x,y)).
$$
If $h$ is a constant function equal to $c>0$, one deduces that 
$$
Ph(x,y)= \frac{c}{|b|} \sharp \{ n\in \nnn(a,b) \ : \ (x,y)\in 
 T_n(\Omega_{n,a})  \}.
$$
One can see that this cardinal number is $|b|$, which proves that there exists a
constant invariant density. 

\noindent
But the theorem proves the existence of an invariant measure  $h^*m$, 
given by $P_1 \1_{\Om} = h^*$. According to Lemma 4.1 of \cite{ITM}, 
$$ P_1 \1_{\Om} =\lim\limits_{n \rightarrow +\infty} \frac{1}{n} \sum\limits_{k=1}^n P^k \1_{\Om} = \1_{\Om}.
$$
Consequently, $h^*$ is a constant function.

%
%

Laboratoire de Math\'ematiques, FR CNRS 3399, EA 4535, Universit\'e de Reims
Champagne-Ardenne, Moulin de la Housse, B. P. 1039, F-51687
Reims, France,

{\it E-mail:} {\tt lisette.jager@univ-reims.fr}

{\it E-mail:} {\tt jules.maes@univ-reims.fr }

{\it E-mail:} {\tt alain.ninet@univ-reims.fr}


\begin{thebibliography}{99}

\bibitem[AFLV]{AFLV}
 ALVES José F.,  FREITAS Jorge M., LUZZATTO Stefano, VAIENTI Sandro,
{\it  From rates of mixing to recurrence times via large deviations,  Advances in Mathematics},  \textbf{228} (2011), n° \textbf{2} 1203-1236.
\bibitem[CE]{CE}
COLLET Pierre, ECKMANN  Jean-Pierre, {\it Concepts and results in chaotic dynamics: a short course}. Theoretical and Mathematical Physics. Springer-Verlag, Berlin (2006).
\bibitem[HK]{HK}
HOFBAUER Franz, KELLER Gerhard, {\it Ergodic properties of invariant measures for piecewise monotonic transformations}, Mathematische Zeitschrift {\bf 180} (1982), 119-140.
\bibitem[ITM]{ITM}
IONESCU TULCEA C.T., MARINESCU G., {\it Théorie ergodique pour des classes d'opérations non complètement continues}, Annals of Mathematics {\bf Vol. 52, n°2} (1950), 140-147.
\bibitem[LM]{LM}
LASOTA Andrzej, MACKEY Michael C., {\it Chaos, fractals and noise : stochastic aspects of dynamics}, Springer Verlag, New York (1998)
\bibitem[LIV]{LIV}
LIVERANI Carlangelo, {\it Multidimensional expanding maps with singularities: a pedestrian approach}, Ergodic Theory and Dynamical Systems {\bf Vol. 33, n°1} (2013), 168-182.
\bibitem[SAU]{SAU}
SAUSSOL Benoît, {\it Absolutely continuous invariant measures for multidimensional expanding maps}, Israel Journal of Mathematics {\bf 116} (2000), 223-248.
\bibitem[TON1]{TON1}
 TONG Howell, {\it Nonlinear time series. A dynamical system approach}.
  With an appendix by K. S. Chan. Oxford Statistical Science Series, \textbf{6}. Oxford Science Publications. The Clarendon Press, Oxford University Press, New York  (1990).

\bibitem[TON2]{TON2}
    TONG Howel,  {\it Nonlinear time series analysis since 1990: some personal reflections}. Acta Math. Appl. Sin. Engl. Ser. 18 (2002), no. \textbf{2}, 177-184.

\bibitem[YOU]{YOU}
YOUNG Lai-Sang, {\it Recurrence times and rates of mixing}, Israel Journal of Mathematics {\bf 110} (1999), 153-188.

\end{thebibliography}
\end{document}